\documentclass[mrdjm]{birkart04}

 \usepackage{a4}
 \usepackage{amsmath}
 \usepackage{amssymb}
\newtheorem{thm}{Theorem}[section]

 \newtheorem{rem}[thm]{Remark}
 \DeclareMathOperator{\Log}{Log}

 \renewcommand{\a}{\alpha}
  \renewcommand{\b}{\beta}
\newcommand{\f}{\Phi}
\title{Integral representation of some functions related to the Gamma
function}
 \author{Christian Berg}
\address{%
Department of Mathematics\\
University of Copenhagen\\
Universitetsparken 5\\
DK-2100 Copenhagen\\
Denmark}
\email{berg@math.ku.dk}
\begin{document}

\begin{abstract}
We prove that the functions  $\f(x)=[\Gamma(x+1)]^{1/x}(1+1/x)^x/x$ and
$\log \f(x)$ are Stieltjes transforms.
\end{abstract}



 \maketitle

\section{Introduction and main results}
\label{sec:intro}

In \cite{F:G} the authors introduce a subclass of the
  completely monotonic functions which they call
   {\it logarithmically completely monotonic}, and the main
 result in \cite{F:G:C} is that the function
\begin{equation}\label{eq:fct}
\f(x)=\frac{[\Gamma(x+1)]^{1/x}}{x}\left(1+\frac{1}{x}\right)^x
\end{equation}
is logarithmically completely monotonic.

We characterize the class of logarithmically completely
 monotonic functions as the infinitely divisible completely monotonic
 functions studied by Horn in \cite{Ho}. We prove that Stieltjes transforms
 (see (\ref{eq:Stieltjes}) below) belong to this class and that  $\f$ and
 $\log \f$ are both
Stieltjes transforms. Each of these  statements imply the result of
 \cite{F:G:C}. The following explicit representations are obtained:
 \begin{equation}\label{eq:St1}
 \log \f(x)=\int_0^\infty\frac{\varphi(s)}{s+x}\,ds,\quad x>0,
 \end{equation}
 where
 \begin{equation}\label{eq:St1a}
      \varphi(s)=\left\{\begin{array}{ll}
      1-s &\mbox{if $0\leq s<1$} \\
      1-n/s & \mbox{if $n\leq s<n+1,n=1,2,\ldots$}
      \end{array}
      \right.
\end{equation}
and
\begin{equation}\label{eq:St1b}
\f(x)=1+\int_0^\infty\frac{h(s)}{s+x}\,ds,\quad x>0,
 \end{equation}
 with
 \begin{equation}\label{eq:St1c}
 h(s)=\frac{1}{\pi}\frac{s^{s-1}}{|1-s|^s|\Gamma(1-s)|^{1/s}}\sin(\pi\varphi(s)),\;s\geq 0.
\end{equation}

Note that the density $\varphi(s)$ takes its values in the interval $[0,1]$,
and this is the clue to the fact that also $\f$ is a Stieltjes transform.
The density $h$ is continuous on $[0,\infty[$ with $h(0)=\exp(-\gamma)$, where
$\gamma$ is Euler's constant, and $h(n)=0$ for $n\in\mathbb N$.

\medskip
Recall that a function $f:\left]0,\infty\right[\rightarrow \mathbb{R}$ is said to
 be completely
monotonic, if $f$ has derivatives of all orders and satisfies
$$
(-1)^n f^{(n)}(x) \geq 0 \quad{\mbox{for all}} \quad{x>0} \quad{\mbox{and}}
\quad{n=0,1,2,...}.
$$
Bernstein's Theorem, cf. \cite[p. 161]{W},
states that $f$ is completely monotonic if and only if
\begin{equation}\label{eq:Bern}
f(x)=\int_{0}^{\infty}e^{-xs}d\mu(s),
\end{equation}
where $\mu$ is a nonnegative measure on $[0,\infty)$ such that the integral
converges for all $x>0$. The set of completely monotonic functions is denoted
$\mathcal C$.

In \cite{F:G:C} the authors call a function $f:]0,\infty[\rightarrow
]0,\infty[$ {\it logarithmically completely monotonic} if it is $C^\infty$ and
\begin{equation}\label{eq:lcm}
(-1)^k[\log f(x)]^{(k)}\geq 0,\;\mbox{for}\; k=1,2,\ldots.
\end{equation}
If we denote the class of logarithmically completely monotonic functions by
$\mathcal L$, we have $f\in \mathcal L$ if and only if $f$ is a positive
$C^{\infty}$-function such that $-(\log f)'\in\mathcal C$.

The functions of class $\mathcal L$ have been implicitly studied in
\cite{A:B2}, and Lemma 2.4(ii) in that paper  can be stated
 as the inclusion
 $\mathcal L\subset\mathcal C$, a fact also established in \cite{F:G}.

The class $\mathcal L$ can be characterized in the
 following way, established by Horn\cite[Theorem 4.4]{Ho}:

 \begin{thm}\label{thm:L} For a function
 $f:]0,\infty[\rightarrow  ]0,\infty[$ the following are equivalent:
 \begin{enumerate}
 \item[(i)] $f\in\mathcal L$
 \item[(ii)] $f^{\a}\in\mathcal C$ for all $\a>0$
 \item[(iii)] $\root{n} \of{f}\in\mathcal C$ for all $n=1,2,\ldots$.
 \end{enumerate}
 \end{thm}

Another way of expressing the conditions of Theorem \ref{thm:L} is that the
functions in $\mathcal L$ are those completely monotonic functions
for which the representing measure $\mu$ in (\ref{eq:Bern}) is infinitely
divisible in the convolution sense: For each $n\in\mathbb N$ there exists a
positive measure $\nu$ on $[0,\infty[$ with $n$'th convolution power
equal to $\mu$, viz. $\nu^{*n}=\mu$. By condition (ii) there exists a
 convolution semigroup $(\mu_{\a})_{\a>0}$ of positive measures such that
 the Laplace transform of $\mu_\a$ is $f^\a$.
 Note that the convolution of any two positive measures on $[0,\infty[$ is
  well-defined and we have $\mu_\a*\mu_\b=\mu_{\a+\b}$.

  In the special case of $f(0+)=1$
 this is very classical: This is the description of infinitely divisible
  distributions in probability. Since there are probabilities which are not
  infinitely divisible we have $\mathcal C\setminus\mathcal L\neq\emptyset$.

In various papers complete monotonicity for special functions has been
established by proving the stronger statement that the function is a Stieltjes
transform, i.e. is of the form
\begin{equation}\label{eq:Stieltjes}
f(x)=a+\int_0^\infty \frac{d\mu(s)}{s+x},
\end{equation}
where $a\geq 0$ and $\mu$ is a nonnegative measure on $[0,\infty[$ satisfying
$$
\int_0^\infty \frac{1}{1+s}d\mu(s)<\infty.
$$ 
See \cite{A:B1},\cite{A:B2},\cite{Be1},\cite{Be2},\cite{B:P1},\cite{B:P2}.

The set of Stieltjes transforms will be denoted $\mathcal S$. We clearly have
 $\mathcal S\subset\mathcal C$. For more information
about this class see \cite{B:F}.

\begin{thm}\label{thm:SinL} $\mathcal S\setminus\{0\}\subset\mathcal L$.
\end{thm}

\begin{thm}\label{thm:main}
The functions
$$
\f(x)=\frac{[\Gamma(x+1)]^{1/x}}{x}\left(1+\frac{1}{x}\right)^x
$$
and
$$
\log \f(x)=\frac{\log\Gamma(x+1)}{x}-\log x+x\log\left(1+\frac{1}{x}\right)
$$
are Stieltjes transforms with the representations (\ref{eq:St1b}) and
 (\ref{eq:St1}).
\end{thm}

\begin{rem} {\rm The class $\mathcal S$ has the following stability
 properties: If $f\in\mathcal S, f\neq 0$ then $1/f(1/x)$ and $1/(xf(x))$
 are again Stieltjes transforms, cf. \cite{Be1}. Therefore the following
 functions belong to $\mathcal S$:
 $$
 \frac{1}{[\Gamma(1+x)]^{1/x}(1+1/x)^x},\quad
 [\Gamma(1+1/x)]^x(1+x)^{1/x},\quad
 \frac{1}{x[\Gamma(1+1/x)]^x(1+x)^{1/x}}.
 $$
 It was proved in \cite{A:B2} that $[\Gamma(1+1/x)]^x\in\mathcal S$, so also
 the following functions are Stieltjes transforms:
$$
\frac{[\Gamma(1+x)]^{1/x}}{x},\quad
\frac{1}{[\Gamma(1+x)]^{1/x}},\quad
\frac{1}{x[\Gamma(1+1/x)]^x}.
$$
In \cite{A:B1} it was proved that $(1+1/x)^{-x}\in\mathcal S$. Therefore
the function $\f$ given by (\ref{eq:fct}) is a quotient of known Stieltjes
transforms, but this does not imply that the function itself is a Stieltjes
transform.}
\end{rem}

\section{Proofs} For completeness we include a proof of Theorem \ref{thm:L}.

\begin{proof}
"(i) $\Rightarrow$ (ii)". Since $f\in\mathcal L$ implies $f^\a\in\mathcal L$
for all $\a>0$, it is enough to prove the inclusion $\mathcal L\subset
\mathcal C$. Although this is done in \cite{A:B2} and \cite{F:G} we include
the easy proof. By assumption $-(\log f)'=-f'/f\in\mathcal C$,
so in particular $-f'\geq 0$. Assume now that $(-1)^kf^{(k)}\geq 0$ for
$k\leq n$. Then
\begin{eqnarray*}
&&(-1)^{n+1}f^{(n+1)}=(-1)^n((-\log f)'f)^{(n)}\\
&=&\sum_{k=0}^n\binom{n}{k}(-1)^k((-\log f)')^{(k)}(-1)^{n-k}f^{(n-k)}\geq 0,
\end{eqnarray*}
and (ii) follows by induction.

\medskip
"(ii) $\Rightarrow$ (iii)" is obvious.

\medskip
"(iii) $\Rightarrow$ (i)". If $f^{1/n}\in\mathcal C$ we have in particular
$-(f^{1/n})'=(-1/n)f^{-1+1/n}f'\in\mathcal C$. Multiplying by $n$ and letting
 $n\to\infty$ we see that the limit function $-f'/f$  belongs to $\mathcal C$,
because $\mathcal C$ is closed under pointwise limits, cf. \cite{B:F}.
 This establishes (i).
\end{proof}

 \medskip
\noindent{\it Proof of Theorem \ref{thm:SinL}:} Let $f\in\mathcal S$ be non-zero, and
 let $\a>0$. By Theorem \ref{thm:L} it is enough to prove that
  $f^\a\in\mathcal C$. Writing $\a=n+a$ with $n=0,1,\ldots$ and $0\leq a<1$
  we have $f^\a=f^nf^a$, and using the stability of $\mathcal C$ under
   multiplication and that $f^a\in\mathcal S$, cf. \cite{Be0}, the assertion
   follows. $\hfill\square$

\medskip
\noindent{\it Proof of Theorem \ref{thm:main}:} Using the
 expression (\ref{eq:St1a}) for $\varphi$ we find
\begin{eqnarray*}
&&\int_0^\infty\frac{\varphi(s)}{s+x}\,ds=\int_0^1\frac{1-s}{s+x}\,ds
+\sum_{k=1}^\infty\int_k^{k+1}\frac{1-k/s}{s+x}\,ds\\
&=&-1+(x+1)\log\left(1+\frac{1}{x}\right)\\
&+&\sum_{k=1}^\infty
\left[\left(1+\frac{k}{x}\right)\log\left(1+\frac{1}{x+k}\right)-\frac{k}{x}
\log\left(1+\frac{1}{k}\right)\right].
\end{eqnarray*}
Therefore,
$$
\int_0^\infty\frac{\varphi(s)}{s+x}\,ds=\log \f(x)
$$
if and only if
\begin{equation}\label{eq:help}
\log\Gamma(x+1)=x(\log(1+x)-1)+\sum_{k=1}^\infty
\left[(k+x)\log\left(1+\frac{1}{x+k}\right)-k\log\left(1+\frac{1}{k}\right)
\right]
\end{equation}
for $x\geq 0$. Both sides vanish for $x=0$, and they have the same derivative
$\psi(x+1)$, where $\psi$ is the digamma function. This follows easily by
the classical formula
$$
\psi(x)=\log x+\sum_{k=0}^\infty\left[\log\left(1+\frac{1}{x+k}\right)-
\frac{1}{x+k}\right]
$$
cf. \cite[8.362(2)]{G:R}.
This shows that $\log \f$ is a Stieltjes transform with the representation
(\ref{eq:St1}). In particular $\f$ is completely monotonic with the
limit 1 at infinity.

To see that a function $f$ is a Stieltjes transform we will use the
 characterization
of these functions via complex analysis, see \cite[p. 127]{Ak} or \cite{Be1}.
It is necessary and sufficient that $f$ has a holomorphic extension to the
cut plane $\mathcal A=\mathbb C\setminus]-\infty,0]$ and satisfies
$\textrm{Im}\,f(z)\leq 0$ for $\textrm{Im}\,z>0$ and $f(x)\geq 0$ for $x>0$.
For a Stieltjes transform $f$ given by (\ref{eq:Stieltjes}) we have
 $a=\lim_{x\to\infty} f(x)$,
and the measure $\mu$ is the limit in the vague topology of
$$
-\frac{1}{\pi}\textrm{Im}\,f(-x+iy)\,dx
$$
as $y\to 0^+$.

We clearly have $\f(x)>0$ for  $x>0$ and
the holomorphic extension of $\f$ 
 is given by $\f(z)=\exp(\log \f(z))$, where $\log \f(z)$ is the holomorphic
extension obtained by the representation (\ref{eq:St1}). This can also be
 described in the following way: For $z\in\mathcal A$ we let
  $\log\Gamma(z)$ be the
  unique holomorphic branch,  which is real for $x>0$, and we let $\Log$
   denote the principal logarithm. Then
\begin{equation}\label{eq:3terms}   
\frac{\log\Gamma(z+1)}{z}-\Log z+z\Log\left(1+\frac{1}{z}\right)
\end{equation}
is a holomorphic branch of $\log \f(z)$ in $\mathcal A$, and since it agrees
with $\log \f(x)$ for $x>0$, we have
 \begin{equation}\label{eq:St2}
 \log \f(z)=\frac{\log\Gamma(z+1)}{z}-\Log z+z\Log\left(1+\frac{1}{z}\right)=
 \int_0^\infty\frac{\varphi(s)}{s+z}\,ds,\quad z\in\mathcal A.
 \end{equation}

 For $z=x+iy,y>0$ we get
 $$
 \textrm{Im}\log \f(x+iy)=-\int_0^\infty \frac{\varphi(s)y}{(s+x)^2+y^2}\,ds,
 $$
 and  since $0\leq\varphi(s)\leq 1$ for $s\geq 0$ we get
$$
\textrm{Im}\log \f(x+iy)\in\left]-\pi,0\right[,
$$
hence
 $$
 \textrm{Im}\,\f(x+iy)=|\f(x+iy)|\sin(\textrm{Im}\log \f(x+iy))<0,
 $$
 which shows that $\f$ is a Stieltjes transform. For $x\geq 0,y\to 0^+$ we
  further get
 \begin{equation}\label{eq:St3}
  -\frac{1}{\pi}\textrm{Im}\,\f(-x+iy)\to h(x):=\frac{1}{\pi}|\f(-x)|
\sin(\pi\varphi(x)),
 \end{equation}
 which is is a continuous nonnegative function on $[0,\infty[$. Therefore the
convergence is uniform for $x$ in compact subsets of $[0,\infty[$, so $h$ is
the density of the representing measure as a Stieltjes transform.

 This shows the following integral representation of $\f$:
 \begin{equation}\label{eq:St4}
     \f(z)=1+\int_0^\infty\frac{h(s)}{s+z}\,ds,\;z\in\mathcal A
 \end{equation}
 where
 $$
 h(s)=\frac{1}{\pi}\frac{s^{s-1}}{|1-s|^s|\Gamma(1-s)|^{1/s}}\sin(\pi\varphi(s)),\;s\geq 0.
 $$
\vskip-1em
\hfill$\square$

\begin{rem} {\rm There is a close relationship between Stieltjes transforms and
Pick functions. For the latter see \cite{Ak} and \cite{D}. It is possible
 to find the integral representation
 (\ref{eq:St1}) of $\log \f(x)$ using integral representations of the three
 terms of (\ref{eq:3terms}). Here $\Log(z)$ and $z\Log(1+1/z)$
 are Pick functions for which the integral  representations are easily
 obtained. The integral representation of the Pick function
  $(\Log\Gamma(z+1))/z$ was found in \cite{B:P2}.
  The author first found the density  (\ref{eq:St1a}) in this way, but
 once $\varphi$ is found the present direct approach seems easier.}
\end{rem}

\end{document}